\documentclass[a4paper, oneside, 11pt]{article}
\RequirePackage[OT1]{fontenc}
\RequirePackage{amsthm,amsmath}
\RequirePackage[numbers]{natbib}
\usepackage{graphics}
\usepackage{latexsym}
\usepackage{amsfonts}
\usepackage{amssymb}
\usepackage[latin1]{inputenc}

\RequirePackage[OT1]{fontenc}
\RequirePackage{amsthm,amsmath}
\RequirePackage[numbers]{natbib}
\RequirePackage[colorlinks,citecolor=blue,urlcolor=blue]{hyperref}
\usepackage{breakurl}
\usepackage{amssymb}
\usepackage[latin1]{inputenc}
\theoremstyle{plain}

\theoremstyle{remark}
\newtheorem{remark}{Remark}[]
\theoremstyle{definition}

\newtheorem{prop}{\bf Proposition}[]

\theoremstyle{remark}
\newtheorem{rema}{Remark}[]
\theoremstyle{definition}

\begin{document}


\title{{\bf{\Large{Exact Good-Turing characterization of the two-parameter Poisson-Dirichlet superpopulation model}}}\footnote{{\it AMS (2000) subject classification}. Primary: 60G58. Secondary: 60G09.}}
\author{\textsc {Annalisa Cerquetti}\footnote{Corresponding author,  E-mail: {\tt annalisa.cerquetti@gmail.com}}\\
\it{\small}\\
  \it{\small }}
\date{\today}
\maketitle{}

\begin{abstract}
Large sample size equivalence between the celebrated {\it approximated} Good-Turing estimator of the probability to discover a species already observed a certain number of times (Good, 1953) and the modern Bayesian nonparametric counterpart has been recently established by virtue of a particular smoothing rule based on the two-parameter Poisson-Dirichlet model. Here  we improve on this result  showing that, for any finite sample size, when the population frequencies are assumed to be selected from a superpopulation with two-parameter Poisson-Dirichlet distribution, then Bayesian nonparametric estimation of the discovery probabilities corresponds to Good-Turing {\it exact} estimation. Moreover under general superpopulation hypothesis the Good-Turing solution admits an interpretation as a modern Bayesian nonparametric estimator under partial information.
\end{abstract}

\section{Introduction}
In 1953 I. J. Good published the results of his collaboration with Alan Turing and their efforts at Bletchley Park to crack German ciphers for the Enigma machine during the World War II. The main result of the celebrated paper {\it The population frequencies of species and the estimation of population parameters} is the {\it approximated} Good-Turing estimator
\begin{equation}
\label{turing}
D_n^{AGT}(l)\approx\frac{l+1}{n}c_{l+1, n}
\end{equation}
of the probability to observe a species already observed $l$ times in a sample of size $n$ assuming {\it virtually nothing} about the underlying population, for $c_{l+1, n}$ the number of species seen $l+1$ times. Good's purpose in writing the paper was to provide a rigorous mathematical setting where to place Alan Turing's first intuition. He eventually shows that (\ref{turing}) is nothing but an approximated version of the {\it exact} Bayesian estimator which could be calculated if the population frequencies were known and that even more accurate approximations may be obtained by suitably smoothing the unknown counting frequencies (cfr Section 8, Good 1953).

In 2016 Favaro {\it et al.} establish an asymptotical equivalence, between the modern Bayesian nonparametric estimator and the celebrated Good-Turing estimator. They demonstrate that, under two-parameter Poisson-Dirichlet $(\alpha, \theta)$ priors (Pitman and Yor, 1997) on the unknown species abundances, the posterior expected probability to discover a species observed $l$ times in a sample of size $n$ agrees, as $n$ tends to infinity, with the {\it approximated} Good-Turing estimator (\ref{turing}) if the frequency counts were smoothed by
\begin{equation}
\label{fav}
c'_{l,n}=\frac{\alpha(1-\alpha)_{l-1}}{l!}k_n,
\end{equation}
for $\alpha \in (0,1)$ and $k_n$ the number of different species observed, i.e.
\begin{equation}
\label{fav2}
D^{BNP}_{n}(l)=D^{AGT}_{n}(l) \approx \frac{l+1}{n} c'_{l+1, n} \approx \frac{l+1}{n}\frac{\alpha(1-\alpha)_{l}k_n}{(l+1)!}=\frac{k_n}{n}\frac{\alpha(1-\alpha)_{l}}{l!}.
\end{equation}

\medskip

The purpose of this paper is to improve on this result by showing that a more in depth analysis of Good's construction allows to prove that, for any {\it finite} sample size, under the same prior hypothesis, the Bayesian nonparametric estimator of the discovery probability corresponds to {\it exact} Good-Turing  estimator (Good, 1953). Moreover, by means of a suitably reformulation of Good's experimental conditions, we point out that, under {\it general} prior superpopulation on the unknown species abundances, the exact Good-Turing  solution admits an interpretation as a modern Bayesian nonparametric solution under {\it partial} information. In other words it turns out that the two estimators just differ for the amount of information from the sample they assume to be available to the experimenter.

\medskip
Let $N=(N_{1,n}, \dots, N_{s,n})$, for $N_{j,n} \geq 0$  be the random abundances vector for a sample of size $n$ taken with replacement from  a multinomial model of fixed parameters (known species relative abundances) $p=(p_1, \dots, p_s)$. Then $N_{j,n} \sim Bin(n, p_j)$ and $P(N_j >0)=1-(1-p_j)^n$. Let $K_n=\sum_{j=1}^s 1\{N_{j,n} >0\}$ be the random number of  different species in the sample and  $(C_{1,n}, \dots, C_{l,n}, \dots, C_{n,n})$ the random vector of the number of species represented $1, \dots, l, \dots, n$ times, then 
\begin{equation}
\label{species}
E_{N}[K_n(p)]=E_{N}[\sum_{j=1}^s 1 \{N_{j,n}> 0\}]=\sum_{j=1}^s [1-(1-p_j)^n]=s-\sum_{j=1}^s (1-p_j)^n,
\end{equation}
and 
\begin{equation}
\label{ewens}
E_{N}[C_{l,n}(p)]=E_{N}[\sum_{j=1}^s 1 \{N_{j,n}=l\}]= \sum_{j=1}^s {n \choose l}p_j^l(1-p_j)^{n-l}.
\end{equation}

Now assume that $P=(P_j)_{j \geq 1}$ is an infinite {\it random} discrete distribution, such that $\sum_ j P_j=1$, let $(\tilde{P}_j)_{j \geq 1}$ be its sized-biased permutation and recall that, for $\tilde{P}_1$ the random relative abundance of the first size-biased pick from $P$ and 
$f$ an arbitrary non-negative measurable function, then
\begin{equation}
E_{P}[\sum_{j=1}^\infty f(P_j)]= E_{\tilde{P}_1}[\tilde{P}_1^{-1} f(\tilde{P}_1)].
\end{equation}
Let $N_1, \dots, N_{K_n}$, for $N_{j} >0$,  be the {\it random} abundances vector of {\it random} size $K_n$ of the different picks in order of appearance observed in a sample of size $n$ taken from $P$ then the analogues of  (\ref{species}) and (\ref{ewens}) yield
\begin{eqnarray}
\label{specieNUMB}
E_{N,P}[K_n(P))]&=& E_P [E_{N|P}[\sum_{j=1}^\infty 1 \{N_{j,n}> 0\}]]\\\nonumber &=& E_{P} [\sum_j [1-(1-P_j)^n]]= E_{\tilde{P}_1}\left[\tilde{P}_1^{-1}[1-(1-\tilde{P}_1)^n]\right]
\end{eqnarray}
and
\begin{eqnarray}
\label{ewensP}
E_{N,P}[C_{l,n}(P)]&=&E_P[E_{N|P}[\sum_{j=1}^\infty 1 \{N_{j,n}=l\}]]\\\nonumber&=& E_{P}[\sum_{j=1}^\infty {n \choose l}P_j^l(1-P_j)^{n-l}]= {n \choose l} E_{\tilde{P}_1}[ \tilde{P}_1^{l-1}(1-\tilde{P}_1)^{n-l}].
\end{eqnarray}

\subsection{Gibbs-type superpopulations}
Gibbs-type models  (Gnedin and Pitman, 2006) are a popular class of {\it species sampling models}  generalizing the random discrete distribution characterizing the random weights of the infinite sum representation of the Dirichlet process (Pitman, 1996, 2003, 2006), while preserving the product form of the induced 
exchangeable partition probability function 
\begin{equation}
\label{kingm}
p(n_1, \dots, n_k)=\sum_{(i_1, \dots, i_k)} {E} \left[  \prod_{j=1}^k P_{i_j}^{n_j} \right],
\end{equation}
where $n_i=|A_i|$,  $(i_1, \dots, i_k)$ ranges over all ordered $k$-tuples of distinct positive integers and $(P_i)_{i \geq 1}$ is any rearrangements of the ranked atoms $(P_i^{\downarrow})_{i \geq 1}$ of $P$. 
Those models are characterized by having (\ref{kingm}) in the product form
\begin{equation}
p_{\alpha, V}(n_1, \dots, n_k)= V_{n,k} \prod_{j=1}^{k} (1 -\alpha)_{n_j-1},
\end{equation}
for $\alpha \in (-\infty, 1)$ and $V=(V_{n,k})$ weights satisfying the backward recursive relation 
\begin{equation}
V_{n,k}=(n -k\alpha)V_{n+1, k} + V_{n+1, k+1},
\end{equation} 
where $V_{1,1}=1$ and $(x)_y=(x)(x+1)\cdots(x+y-1)$ is the usual notation for rising factorials. 
Given $(n_1, \dots, n_k)$ the multiplicities of the first $k$ species observed in a random sample of size $n=\sum_i n_i$,  the probabilities to observe the $j$th {old} species or a {new} yet unobserved species at step $n+1$ are well-known to correspond respectively to $p_{\alpha, V}(n_j^{+})= {V_{n+1,k}} (n_j -\alpha)/{V_{n,k}}$ and $p_{\alpha, V}(n^{k+1})=V_{n+1, k+1}/V_{n,k}$.

The $r$-th falling factorial moment of $C_{l,n}$ under a general Gibbs-type model $(\alpha, V)$ has been derived in Cerquetti (2013) (see also Favaro {\it et al.} 2013) and is given by
\begin{equation}
\mathbb{E}_{\alpha, V}\left[(C_{l,n})_{[r]}\right]= \frac{n! [(1-\alpha )_{l-1}]^r}{(l!)^r  (n-lr)!} \sum_{k-r=0}^{n-rl}
 V_{n,k} S_{n-rl, k-r}^{-1, -\alpha},
\end{equation}
where $S_{n,k}^{-1, -\alpha}$ stands for the  generalized Stirling number of parameters $(-1, -\alpha)$.
For $r=1$, 
\begin{equation}
\label{GibbsCL}
E_{\alpha, V}[C_{l,n}(P)]= {n \choose l} (1-\alpha)_{l-1} \sum_{k-1=0}^{n-l} V_{n,k}S_{n-l, k-1}^{-1, -\alpha}
\end{equation}
and the {\it prior expected number of different species} in a sample of size $n$ specializes as 
\begin{eqnarray}
\label{nublock}
E_{\alpha, V}[K_n(P)]&=&\sum_{l=1}^n E_{\alpha, V}[C_{l,n}(P)]\\\nonumber &=& \sum_{l=1}^{n} {n \choose l} (1-\alpha)_{l-1}\sum_{k-1=0}^{n-l}V_{n,k} S_{n-l,k-1}^{-1, -\alpha}.\\\nonumber
\end{eqnarray}
Notice that both (\ref{GibbsCL}) and (\ref{nublock}) may be derived directly from (\ref{specieNUMB}) and (\ref{ewensP}) if the law  of $\tilde{P}_1$, the {\it structural distribution} of the prior $P=(P_i)_{\geq 1}$ and its expected value are available in closed form.  

\section{Reformulating exact Good-Turing estimation}

Exact Good-Turing estimation of the probability to discover a species encountered $l$ times in a sample of size $n$ is introduced in Good (1953) as follows:
\begin{quote}{\it Suppose that a species is selected at random or rather equiprobably from an infinite population of $s$ species with fixed unknown population frequencies $p_1, \dots, p_s$ and that then its number of occurrences in a sample of size $n$ is $l$. Let $H$ be the statistical hypothesis on $(p_1, \dots, p_s)$. Thus the initial ({\it prior}) probability that the species is the $j$-th of population frequency $p_j$ is $1/s$ and the likelihood that the observed number of occurrences is $l$ is ${n \choose l} p_j^l(1-p_j)^{n-l}$. 
}
\end{quote}

Notice that setting a uniform prior $1/s$ for the probability of every single species and then conditioning to the event that the species $j$ has been observed $l$ times corresponds to setting a prior $(p_1, \dots, p_s)$ of the probability of the single species and conditioning to $l-1$ observations of the same species $j$ observed in the first sampled individual. It follows that Good's experimental setting can be equivalently reformulated as follows: \\

Assume to have a infinite population $p$ of $s$ species, with fixed unknown population frequencies $p_1, \dots, p_s$. Without loss of generality let $\tilde{P}_1$ be the random population frequency of the {\it first} encountered species, taking values $(p_1, \dots, p_s)$ with {\it prior} probability respectively $(p_1, \dots, p_s)$. Now, consider the random variable $N_{1,n}-1$, i.e. the number of the observations of the {\it first} observed species, in a sample of additional $n-1$ observations. Then  $N_{1,n}-1| \tilde{P}_1=p_\mu \sim Bin (n-1, p_\mu)$ and the likelihood that the total number of occurrences of the first species in a random sample of size $n$ is $l$ is given by ${n -1 \choose l-1} p_{\mu}^{l-1} (1 - p_{\mu})^{n-l}$. 

Then Good writes:
\begin{quote}
{\it We shall consider the {\it  posterior} expectation of $Q_l$, the unknown population frequencies of an arbitrary species that is represented $l$ times in the sample.}
\end{quote}
We can then rewrite Good's exact result (cfr. eq. 12 in Good, 1953) for the {\it posterior} probability of the population frequency of a species encountered $l$ times in a sample of size $n$ as follows
$$
Pr(Q_l=p_\mu| H) = Pr \{\tilde{P}_{1}={p}_\mu| N_{1,n}=l, H\}= \frac{p_\mu^l (1-p_\mu)^{n-l}}{\sum_{j=1}^s p_j^{l}(1-p_j)^{n-l}}.
$$
By an application of (\ref{ewens}) the posterior expectation (cfr eq. 15 in Good, 1953) results

\begin{equation}
\label{good1}
E_{p}[\tilde{P}_{1}|N_{1,n}={l}, H]= \frac{\sum_{j=1}^s p_j^{l+1} (1-p_j)^{n-l}}{\sum_{j=1}^s p_j^{l}(1-p_j)^{n-l}}= \frac{(l+1)}{(n+1)}\frac{E_{N+1}(C_{l+1, n+1}(p)|H)}{E_{N}(C_{l,n}(p)|H)},
\end{equation}
for $C_{l,n }$ the random number of species represented $l$ times in a sample of size $n$.

\medskip

Assuming  {\it virtually nothing} about the underlying population $H(p_1, \dots, p_s)$, hence being impossible to evaluate (\ref{good1}), Good derives the estimator originally devised by Alan Turing observing that it makes sense to approximate (\ref{good1}) with  
\begin{equation}
\label{tur}
E_{p}[(\tilde{P}_{1})|N_{1,n}={l}] \approx \frac{l+1}{n}\frac{c_{l+1, n}}{c_{l, n}} 
\end{equation}
for $c_{l,n}$ the number of species actually observed $l$ times in the $n$ sample. Then discusses the possibility to obtain even more accurate approximations by suitably smoothing the $c_{l,n}$'s (cfr Section 3, Good 1953).

\rema{It is worth to notice  at that point that in (\ref{good1}) Good is {\it just} conditioning  on the observed frequency of the single species under consideration, nevertheless in (\ref{tur}) he assumes the whole vector of sample frequencies to be known, being impossible otherwise for $c_{l,n}$ and $c_{l+1,n}$} to be available. We will see in what follows that this discrepancy about the information available to the experimenter is what differentiates Good's approach from the modern bayesian approach and how it vanishes under the two-parameter Poisson-Dirichlet superpopulation hypothesis.

\section{Main results}

Few lines later his main result Good (1953) clarifies:
\begin{quote}
{\it If a specific (statistical) hypothesis $H(p_1, \dots, p_s)$ is accepted it is clearly not necessary to use the approximations since equation (\ref{good1}) can then be used directly. Similarly, if $H$ is assumed to be selected from a superpopulation, with an assigned probability density, then again it is theoretically possible to dispense with the approximation}.
\end{quote}
Let suppose that $H(p)$ is assumed to be selected from a superpopulation $H(P)$ for $P=(P_j)_{j \geq 1}$  a {\it random} discrete distribution model for the unknown population frequencies, that here, without loss of generality, we assume to be countable infinite. Then an application of (\ref{good1}) yields 

\begin{eqnarray}
\label{goodBNP2}
E _P[\tilde{P}_1|N_{1,n}=l]&=&\frac{E_{P}[[\sum_{j=1}^\infty P_j^{l+1}(1-P_j)^{n-l}]}{E_{P}[[\sum_{j=1}^\infty P_j^l(1-P_j)^{n-l}]}\\
&=&\frac{E_{\tilde{P}_1} [\tilde{P}_1^l (1 -\tilde{P}_1)^{n-l}]}{E_{\tilde{P}_1} [\tilde{P}_1^{l-1} (1 -\tilde{P}_1)^{n-l}]}\\\nonumber\ &=&\frac{l+1}{n+1}\frac{E_P[E_{N|P} (C_{l+1, n+1} (P))]}{E_{P}[ E_{N|P} (C_{l,n}(P))]}
\end{eqnarray}
where the last equivalence follows from (\ref{ewensP}). An application of (\ref{goodBNP2}) and (\ref{GibbsCL}) provides the following result

\prop{\it 
For $P$ a general Gibbs-type superpopulation of parameters $(\alpha, V)$ then the exact Good-Turing estimator of the probability to encounter a species already observed  $l$ times in a sample of size $n$ is given by
\begin{equation}
\label{goodBAY_gibbs}
E_{\alpha, V}[\tilde{P}_1| N_{1,n}=l]= (l -\alpha)\frac{\sum_{k-1=0}^{n-l} V_{n+1,k} S_{n-l, k-1}^{-1, -\alpha}}{\sum_{k-1=0}^{n-l} V_{n,k} S_{n-l, k-1}^{-1, -\alpha}}.
\end{equation}
}

\remark{It follows from (\ref{goodBAY_gibbs}) that under Gibbs-type superpopulations the {\it exact} Good-Turing  estimator turns out to be a modern Bayesian nonparametric estimator under {\it partial} information. Notice in fact that the fundamental difference between (\ref{goodBAY_gibbs}) and the Bayesian nonparametric estimator which follows from the sequential rule defining the Gibbs EPPF
\begin{equation}
\label{bnpl}
E_{\alpha, V}(\tilde{P}_j|n_1, \dots,n_j=l,\dots, n_k)= (l -\alpha)\frac{V_{n+1, k}}{V_{n,k}} 
\end{equation}
relies on the conditioning argument. While Good is just conditioning on the information on the observed frequency of the species he wants to predict $(n_1=l)$, the modern Bayesian assumes to know the whole sequence of the multiplicities of the encountered species $(n_1, \dots, n_k)$ and therefore $k_n=\sum_i c_{i,n}$ the number of different species seen. This explains why in (\ref{goodBAY_gibbs}) we are marginalizing both the numerator and the denominator over the number $k_n$ of different species seen in a sample of size $n$ when at least one species has been observed $l$ times }

\medskip

In the large $(\alpha, V)$ Gibbs-type family of superpopulations the two-parameter $(\alpha, \theta)$ Poisson-Dirichlet model is  well-known to be characterized by the peculiar EPPF
\begin{equation}
\label{EPPF_PD}
p(n_1, \dots, n_k)= \frac{(\theta +\alpha)_{k-1 \uparrow \alpha}}{(\theta +1)_{n-1}} \prod_{j=1}^k (1 -\alpha)_{n_j -1},
\end{equation}
for $\alpha \in (0,1)$ and $\theta > -\alpha$ or $\alpha <0$ and $\theta= k|\alpha|$, $k=1, 2, 3,\dots$. It is easy to see that the characteristic form  $g(k)/h(n)$ of the $PD(\alpha, \theta)$ Gibbs weights characterizes the Poisson-Dirichlet family as the only Gibbs-type model for which  Good-Turing exact estimation exactly agrees with the corresponding {\it modern} Bayesian nonparametric estimator.

\prop{\it 
Let $P$  be distributed according to a two-parameter Poisson-Dirichlet  $(\alpha, \theta)$ superpopulation model on the infinite dimensional simplex,  with  structural distribution  $\tilde{P}_1 \sim beta (1-\alpha, \theta +\alpha)$ for $\alpha \in [0,1)$ and $\theta > -\alpha$ or $\alpha <0$ and $\theta= |\alpha|s$ on the $s$-dimensional simplex. Then, for any finite sample size $n$ the exact Good-Turing estimator agrees with the modern Bayesian nonparametric estimator.
}

\proof{

Specializing (\ref{bnpl}) for the $PD(\alpha, \theta)$ model
\begin{equation}
\label{ellePD}
E_{\alpha, \theta}(\tilde{P}_j|n_1, \dots, n_j=l, \dots, n_k)= (l -\alpha)\frac{V_{n+1, k}}{V_{n,k}} = \frac{(\theta + \alpha)_{k-1 \uparrow \alpha}}{(\theta +1)_{n}}  \frac{(\theta +1)_{n-1}}{(\theta + \alpha)_{k-1 \uparrow \alpha}}= \frac{l-\alpha}{(\theta +n)}.
\end{equation}
By elementary results for the beta distribution, specializing (\ref{goodBAY_gibbs}) yields

\begin{equation}
\label{pdat_good}
E_{\alpha, \theta} [\tilde{P}_1|N_{1,n}=l]=\frac{(1 -\alpha)_{l} (\theta +\alpha)_{n-l}}{(\theta +1)_n} \frac{(\theta +1)_{n-1}}{(1-\alpha)_{l-1}(\theta +1)_{n-l}}= \frac{(l-\alpha)}{(\theta +n)}.
\end{equation} 
}

\begin{rema}
For $\alpha < 0$ and $\theta= |\alpha|s$ then (\ref{pdat_good}) yields Johnson (1932) estimate: $(l+|\alpha|)/(n +|\alpha|s)$ under the assumption that $P_1, \dots, P_s$ is selected from the finite dimensional simplex with {\it symmetric Dirichlet} distribution of parameters $(|\alpha|, \dots, |\alpha|)$.   For $\alpha=-1$ then (\ref{pdat_good}) yields Jeffreys (1948) estimate under {\it uniform} distribution on the $s$-dimensional simplex. Notice that both results are cited in Good (1953) (cfr. pag. 241) as examples of assuming $H$ to be selected from a superpopulation. 
Notice that Good underlines that both Jeffreys' and Johnson's solutions ignore any information that can be obtained from the entire set of frequencies of all species, but this is just what characterized the {\it exact} estimator he proposes. He highlights the importance of the whole frequencies vector for the smoothing of the $c_{l,n}$ when one needs to improve the performances of the {\it approximated} estimator, nevertheless, when building the exact solution he do not assume this information to be available. With no surprise the two-parameter Poisson-Dirichlet superpopulation model with $\alpha \in (0,1)$ corresponds to model $H_1$ in Good (1953) Section 7 (pag. 248) for $\alpha^*= -\alpha-1$ and $\beta= \theta+\alpha-1$.
\end{rema}

\section*{References}
\newcommand{\bibu}{\item \hskip-0.4cm}
\begin{list}{\ }{\setlength\leftmargin{0.5cm}}


\bibu \textsc{Cerquetti, A.} (2013) Marginals of multivariate Gibbs distributions with application in Bayesian species sampling. {\it Electronic Journal of Statistics, 7, 697--716}
\bibu \textsc{Favaro, S. Lijoi, A. and Pr\"unster, I.} (2013) Conditional formulae for Gibbs-type exchangeable random partitions. {\it Ann. Appl. Probab.}, {\bf 23}, 1721--1754
\bibu \textsc{Favaro, S., Nipoti, B., Teh, Y.W.} (2016) Rediscovery of Good-Turing estimators via Bayesian nonparametrics. {\it Biometrics}, 72, 136-145. 

\bibu \textsc{Good, I. J.} (1953) The population frequencies of species and the estimation of population parameters. {\it Biometrika}, 40, 237--264

\bibu \textsc{Gnedin, A. and Pitman, J. } (2006) {Exchangeable Gibbs partitions  and Stirling triangles.} {\it Journal of Mathematical Sciences}, {\bf 138}, 3, 5674--5685.

\bibu \textsc{Jeffreys, H.} (1948) {\it Theory of probability}, 2nd Ed. Oxford: Clarendon Press. 
\bibu \textsc{Johnson, W.E. } (1932) Appendix (edited by R.B. Braithwaite) to {\it Probability: deductive and inductive problems}, {\it  Mind}, {\bf 41}, 421-423.

\bibu \textsc{Pitman, J.} (1996) Some developments of the Blackwell-MacQueen urn scheme. In T.S. Ferguson, Shapley L.S., and MacQueen J.B., editors, {\it Statistics, Probability and Game Theory}, volume 30 of {\it IMS Lecture Notes-Monograph Series}, pages 245--267. Institute of Mathematical Statistics, Hayward, CA.

\bibu \textsc{Pitman, J.} (2003) {Poisson-Kingman partitions}. In D.R. Goldstein, editor, {\it Science and Statistics: A Festschrift for Terry Speed}, volume 40 of Lecture Notes-Monograph Series, pages 1--34. Institute of Mathematical Statistics, Hayward, California.

\bibu \textsc{Pitman, J.} (2006) {\it Combinatorial Stochastic Processes}. Ecole d'Et\'e de Probabilit\'e de Saint-Flour XXXII - 2002. Lecture Notes in Mathematics N. 1875, Springer.

\bibu \textsc{Pitman, J. and Yor, M.} (1997) The two-parameter Poisson-Dirichlet distribution derived from a stable subordinator. {\it Ann. Probab.}, {\bf 25}, 855--900.

\end{list}
\end{document}